\documentclass[a4paper,12pt]{article}

\usepackage{geometry}

\usepackage{amsmath}
\usepackage{amssymb}
\usepackage{amsthm}
\usepackage[latin1]{inputenc}
\newtheorem{theorem}{Theorem}

\newtheorem{lemma}{Lemma}

\newtheorem{definition}{Definition}

\newcommand{\re}{\mathbb{R}}
\newcommand{\nat}{\mathbb{N}}

\newcommand{\rar}{\rightarrow}
\newcommand{\stone}{\mathcal{M}_T (X)}
\newcommand{\sttwo}{\mathcal{M}_S (Y)}
\newcommand{\lip}{\mathcal{L}(X,\mathbb{R})}
\newcommand{\lipy}{\mathcal{L}(Y,\mathbb{R})}
\def\cqd {\,  \begin{footnotesize}$\square$\end{footnotesize}}

\begin{document}

%\title{Ergodic Optimization, Thermodynamic Formalism and Ergodic Games}
\title{Ergodic and Thermodynamic Games}
\author{Rafael Rigão Souza} % Autor
%\thanks{The author is partially supported by FAPERGS (proc. 002063-2551/13-0).}
\date{\today}

\maketitle

\begin{abstract}

Let $X$ and $Y$ be compact sets, and $T:X\to X $, $S:Y \to Y$ be continuous maps. Let $$\varphi_i(\mu,\nu)=\int_{X \times Y} A_i(x,y) d\mu(x) d\nu(y)\;\;\mbox{for} \;\; i=1,2,$$ where $\mu$ is $T$-invariant and $\nu$ is $S$-invariant, be {\it payoff functions} for a game (in the usual sense of game theory) between players that have the set of invariant measures for $T$ (player 1) and $S$ (player 2) as possible {\it strategies}. Our goal here is to establish the  notion of Nash equilibrium for the game defined by these payoffs and strategies. The main tools come from ergodic optimization (as  we are optimizing over the set of invariant measures) and thermodynamic formalism (when we add to the integrals above the entropy of measures in order to define a second case to be explored). Both cases are ergodic versions of non-cooperative games. 
We show the existence of Nash equilibrium points with two independent arguments. 
One of the arguments deals with the case with entropy, and uses only tools of thermodynamical formalism, while the other, that works in the case without entropy but can
be adapted to deal with both cases, uses the Kakutani fixed point. We also present examples and briefly discuss uniqueness (or lack of uniqueness). In
the end, we present a different example where players are allowed to collaborate. This final example shows connections between cooperative games and
ergodic transport.

%One of them, for the case with entropy, using only tools of thermodynamical formalism and the other, that works in both cases, using the Kakutani fixed point along with tools of ergodic optimization. We also present examples and prove uniqueness in the case with entropy (and discuss uniqueness or lack of uniqueness in the general game theory framework). In the end we present a different example where players are allowed to collaborate. In this final case we use tools from ergodic transport. 

\end{abstract}

{\footnotesize The author is partially supported by FAPERGS (proc.002063-2551/13-0).}

%\thispagestyle{empty}

 % \newpage

\section{Introduction}

%We consider here a non-cooperative {\it ergodic game}, where the usual ingredients that define a game are given in terms of ergodic concepts, making use of ideas that came from ergodic optimization and thermodynamic formalism.

%OPCAO 2: 
A strategic non-cooperative game (in the usual sense of game theory - see \cite{NM}) is defined by a set of players, and for each player: (i) a set of strategies and (ii) a payoff (or profit) function that depends on the strategies chosen by each one of the players. Here, we consider ergodic versions of such games.  Suppose there are two continuous maps $T:X \to X$ and $S:Y\to Y$ defined on compact metric spaces $X$ and $Y$.
%We suppose player 1 has as his set of strategies the set of invariant measures for $T$, while the second player has as strategies the set of invariant measures for $S$. 
Suppose  we have two players,  and the strategies of player 1 are given by the set of invariant measures for $T$, while the strategies of player 2 are the invariant measures for  $S$. 
If $A_i:X \times Y \to \mathbb{R} \;,\; i=1,2$, are two Lipschitz potentials, we  define, in a first model,  the {\it payoff function} for player $i$ as the function
$$
\varphi_i(\mu,\nu)=\int_{X \times Y} A_i(x,y) d\mu(x) d\nu(y)\;\;\mbox{for} \;\; i=1,2,
$$ 
where $\mu$ is $T$-invariant and $\nu$ is $S$-invariant. In a second model, we  define as payoffs
\begin{equation*}
\begin{cases}
\varphi_1(\mu,\nu)=\int_{X \times Y} A_1(x,y) d\mu(x) d\nu(y)+h(\mu),
\\
\varphi_2(\mu,\nu)=\int_{X \times Y} A_2(x,y) d\mu(x) d\nu(y)+h(\nu),
\end{cases}
\end{equation*}
where  $h(\mu)$ and $h(\nu)$ are the metric entropies of $\mu$ and $\nu$. 
The first case, which is the ergodic optimization  setting, 
 will define what we call an ergodic game, while the second case, which is related to the thermodynamical formalism, will define what can be called an ergodic game (as in the first case) or a thermodynamic game. % Nevertheless, we will call both kind of games as ergodic games.

%In both cases we consider the product of measures when integrating potentials $A_i$, which reflects the lack of collaboration between players. Therefore we are dealing with non-cooperative games defined by ergodic or dynamic objects. The concept of 
%The first model is in the ergodic optimization  setting while the second is related to the thermodynamical formalism.

%\medskip

% We can also consider the entropy of measures added to the integral of potentials when defining the payoff functions. If we do not consider entropy we have an application of ergodic optimization (case 1) and if we also have entropy we have an application of thermodynamical formalism (case 2)

%The game we are considering here is the usual strategic game in normal form between a finite number of players (see \cite{NM}). 

We consider two players, but the general case of $N$ players is an easy generalization. The ergodic concepts are used to define the space of strategies for each player.
In both cases, we consider the product of measures when integrating potentials $A_i$, which reflects the lack of collaboration between players: we suppose, as usual in the non-cooperative game theory, that players do not communicate with each other. We present a concept of equilibrium (or solution) for the game based on the Nash equilibrium used in non-cooperative game theory, where players are in equilibrium if any one of them who change his strategy gets a worst result provided the other players keep their strategies unchanged (which shows that a Nash equilibrium point is stable, in some sense). 

We present two independent proofs of existence of Nash equilibrium. The first one is more general and just requires the potentials to be Lipschitz continuous,  can be applied to both cases (with or without entropy), and uses a set-valued fixed point due to Kakutani, which is not very usual in the ergodic theory literature. The second proof uses more conventional arguments of thermodynamical formalism but can only be applied to the second case (payoffs with entropy - or thermodynamic games). 

Both proofs are  independent and perhaps the reader more familiar to ergodic optimization will be more comfortable with the second proof, while the first one is more general but uses the Kakutani theorem, which was used in game theory to provide existence of equilibrium of classical (non-ergodic) games, but, as already said, is not very usual in the areas of ergodic optimization and thermodynamic formalism.
%We also present a proof of uniqueness of Nash equilibrium in the case with entropy, using arguments of thermodynamical formalism. 
%The Wasserstein-1 metric plays an important role in both existence and uniqueness proofs in this case.

The Wasserstein-1 metric plays an important role in the proofs of existence.

We present some examples in the text and also discuss why uniqueness of solutions is not a trivial question.
We finish with a different model (the hierarchical case), that is a special case without entropy, that allow us to use some concepts of cooperative games leading to the use of ergodic transport techniques (see \cite{LM,LMMS2}). We stress, however, the fact that this example is unique in the text: all other situations considered here, including the definition of Nash equilibrium and proofs of existence, are in the domain of non-cooperative {\it ergodic} game theory, where players are supposed not  to communicate or collaborate with each other. % We do not elaborate much in this last example: we content ourselves in presenting the example, which can be deeper analysed in a future work.

%\medskip

%\footnote{Bibliographical Revision...}

This paper is structured as follows: in the second section of the paper, we establish the notion of ergodic games %after a brief introduction to 
in the two cases we consider here, % (without and with entropy - the first can be called an ergodic game while the second can be called a thermodynamic game - but the terminology ergodic game could be also used in the second case), 
and present some examples.
In the third section, we prove existence of Nash equilibrium for 
the first case (ergodic games, without entropy), while in the fourth section, we prove existence for the second case (thermodynamic games). In the fifth section, we discuss uniqueness. In the last section,  we present a final example that uses ideas of ergodic transport.

\medskip

%\footnote{ OPCAO 1 para PRIMEIRO paragrafo: Here we will use ideas of ergodic optimization and thermodynamical formalism to analyse a game (in the usual sense of game theory)  between two players that have as  strategies invariant measures for a map $T:X \to X$ (player 1) or invariant measures for a map $S:Y\to Y$ (player 2). We suppose both maps $T:X \to X$ and $S:Y\to Y$ are continuous and both sets $X$ and $Y$ are compact metric spaces. As payoffs (or profit functions) we consider two cases: one that uses just a potential $A_i:X \times Y \to \mathbb{R}$ to define the payoff of the $i$-th player (by integrating this payoff with respect to the product measure of both strategies choosen by the players) and other case that adds to the integral of the potential $A_i$  the metric entropy of the strategy of player $i$. The first model is in the ergodic optimization  setting while the second is related to the thermodynamical formalism.}

%\footnote{Maybe an additional section with some uniqueness results (if time allows...)}

%\subsection{Games}

\section{Ergodic theory and games} 
 
\medskip
%we very briefly introduce the basic facts of non-cooperative game theory. 

Before setting up the ergodic games, we recall that a strategic non-cooperative game in normal form (see \cite{NM}) between a finite number of players is given by a set  $\{1,2,...,N\}$ (set of players), and for each player $i$: 

(a) a set $S_i$ , called the set of {\it strategies} for player $i$;

(b) a function $\varphi_i:S_1 \times ... \times S_N \to \mathbb{R}$, called the {\it payoff} or {\it utility function} of player $i$.

Given such a game, each player wants to maximize its own payoff. We suppose that payoff functions depend on the choices of strategies of {\it all} players, which explains the interaction scheme we have here. We also suppose players do not communicate / collaborate with each other. 
Therefore, when we search for the best strategy for each player, the situation is totally different from the case where one just has to optimize a certain function. 
There are different concepts of solution or equilibrium points for games (dominated strategies, Pareto optimality, etc). We will use here one of the most important concepts of equilibrium for non-cooperative games, due to Nash (see \cite{Na}), which will be carefully introduced after we set up the ergodic (and thermodynamic) games.

  Let $(X,d_X)$ and $(Y,d_Y)$ be two compact metric sets, and  $T:X\rar X$ and $S:Y \rar Y$ be two continuous maps. Let $\mathcal{M}(X)$ and $\mathcal{M}(Y)$ represent respectively  the  probability measures on the Borel sets of $X$ and $Y$, while
 $\mathcal{M}_T(X)$ and $\mathcal{M}_S(Y)$ represent, respectively,  the invariant probability measures for $T$ and $S$.

The following lemma is a well-known result (see \cite{Ka}): 
\begin{lemma}\label{primeiro}
The sets  $\stone$ and $\sttwo$ are, respectively, non-empty compact and convex subsets of 
$\mathcal{M}(X)$ and $\mathcal{M}(Y)$.
\end{lemma}

Here, as usual, we are using the weak-* topology.
Therefore,  $\mu_n \in \mathcal{M}(X)$ converges to $\mu \in \mathcal{M}(X)$ if and only if $\int \psi d\mu_n \to \int \psi d\mu$ for any continuous function $\psi:X\to \mathbb{R}$, and we have the analogous definition  for sequences of measures in $\mathcal{M}(Y)$.

We  define an  {\it ergodic game} between two players 1 and 2 by considering {\it payoffs} 

\[
\begin{cases}
\varphi_1:\mathcal{M}_T(X)\times \mathcal{M}_S(Y)\to \mathbb{R},\\ 
\varphi_2:\mathcal{M}_T(X)\times \mathcal{M}_S(Y)\to \mathbb{R}, 
\end{cases}
\]
%$$\varphi_1:\mathcal{M}_T(X)\times \mathcal{M}_S(Y)\to \mathbb{R}$$ and
%$$\varphi_2:\mathcal{M}_T(X)\times \mathcal{M}_S(Y)\to \mathbb{R},$$
where  $\mathcal{M}_T(X)$ can be seen as the set of strategies for player 1, while  $\mathcal{M}_S(Y)$ is the set of strategies for player 2.

In this paper, we will consider that the payoff functions are given by 
\begin{equation}\label{potenciais}
\varphi_i(\mu,\nu)=\int_{X \times Y} A_i(x,y) d\mu(x) d\nu(y)
\,,\;\;i=1,2,
\end{equation}
or %the payoffs are given by
\begin{equation}\label{payoff-entropia}
\begin{cases}
\varphi_1(\mu,\nu)=\int_{X \times Y} A_1(x,y) d\mu(x) d\nu(y)+h(\mu),
\\
\varphi_2(\mu,\nu)=\int_{X \times Y} A_2(x,y) d\mu(x) d\nu(y)+h(\nu),
\end{cases}
\end{equation}
where in both cases $A_i:X \times Y \rar \mathbb R$, for each $i=1,2$ are  Lipschitz continuous {\it potentials}, and in the second case  $h(\mu)$ and $h(\nu)$ are the metric entropies of $\mu$ and $\nu$. 
In  the product space $X \times Y$, we use the metric given by  $d((x,y)),(x',y'))=d_X(x,x')+d_Y(y,y')$.
We call both  games  ergodic games, and the second can also be called a thermodynamic game.

It is important to remark that the only measures we consider here are the product of invariant measures. This structure that only allows product measures is necessary because of the lack of cooperation between players. In the last section,  we will present an entirely different game where players can cooperate, and there we will use measures in $X \times Y$ that are no longer product measures.

% We also remark that we will study  here very general potentials.  We will just require  Lipschitz continuity for the potentials. 
%The usual concept of equilibrium in the general case is the one due to Nash. %And we will see that defining Nash equilibrium involves two simultaneous ergodic optimization problems.

In order to introduce Nash equilibrium in ergodic games,
we begin by studying the best response of player $1$ given that player $2$ is using a strategy $\nu\in \sttwo$. Here,
in the case payoffs are given by \eqref{potenciais},
the search for a best response set will involve a usual ergodic optimization problem, while in the case that 
payoffs are given by \eqref{payoff-entropia}, we will search for equilibrium measures. 
 However, a solution for the game involves the solution of  two simultaneous 
ergodic optimization problems (or the search for two equilibrium measures), %as we have two players, and each problem depends on the other, 
in a sense that soon will be clear. 
Note that, in the case 
payoffs are given by \eqref{potenciais}, i.e., in the case
$$\varphi_i(\mu,\nu)=\int_{X \times Y} A_i(x,y) d\mu(x) d\nu(y)=$$
$$=\int_X \left(\int_Y  A_i(x,y)  d\nu(y) \right) d\mu(x)\,,\;\;i=1,2,$$ 
we see that, for $\nu$ fixed, the map $\mu \to \varphi_1(\mu,\nu)$ is a linear functional and, therefore, there exist at least one measure $\mu$ that maximizes 
\begin{equation}\label{funcional}
 \mu \mapsto \varphi_1(\mu,\nu).
\end{equation}
 So, for any  $\nu \in \sttwo$, we define  $BR_1(\nu)$ as the set of maximum points for \eqref{funcional}, i.e.,
$$BR_1(\nu)=\left\{\mu \in \mathcal{M}_T(X) \;\|\; \varphi_1(\mu,\nu)\geq \varphi_1(\eta,\nu) \; \forall \eta \in \mathcal{M}_T(X) \right\},$$
  which is a compact and convex subset of $\stone$.
 $BR$ stands for {\it best response}, a terminology used in game theory, and $BR_1(\nu)$ can be seen as the set of optimal strategies for player $1$ supposing player $2$ is using the strategy $\nu$.
%IMPORTANTE IMPORTANTE IMPORTANTE Olhar se no fundo, quando fixamos $\nu$ e obtemos a melhor estratégia $\mu$, não temos um problema de transporte usaul. Teria que ver o problema de transporte nos dois sentidos. O que a slackness condition tem a ver com isso tudo ???
Analogously, we define $BR_2(\mu)$ using the second potential $A_2$: if $\mu\in \stone$,
$$BR_2(\mu)=\left\{\nu \in \mathcal{M}_S(Y) \;\|\; \varphi_2(\mu,\nu)\geq \varphi_2(\mu,\eta,) \; \forall \eta \in \mathcal{M}_S(Y) \right\},$$ 
  which is a compact and convex subset of $\sttwo$.
  
We remark the fact that $BR_1(\nu)$ is a subset of $\stone$, while 
$BR_2(\mu)$ is a subset of $\sttwo$. In fact, if $M$ is a set and $\mathcal{P}(M)$ denotes the set of subsets of $M$ (sometimes denoted in the literature by $2^M$), we can write
\[
\begin{cases}
BR_1:\sttwo \to \mathcal{P}(\stone), \\ 
BR_2:\stone \to \mathcal{P}(\sttwo).
\end{cases}
\]

%The proof of the following lemma is straightforward and will be ommited:
We summarize these results in 
\begin{lemma}\label{BR_properties}
In the case 
payoffs are given by \eqref{potenciais},
\begin{enumerate}
\item [a)] For any $\nu$ in 
$\sttwo$, $BR_1(\nu)$ is a convex and non-empty compact subset of $\stone$.
\item[b)] For any $\mu$ in 
$\stone$, $BR_2(\mu)$ is a convex and non-empty compact subset of $\sttwo$.
\end{enumerate}
\end{lemma}

In the case 
payoffs are given by \eqref{payoff-entropia}, the best-response sets can analogously be defined, and a version of the lemma above also holds under some additional hypothesis, see lemma \ref{TFG} in section \ref{sec:entropia}.

Now we are ready to define Nash Equilibrium. The following definition works in both cases of potentials we are considering here. We  remark that Nash Equilibrium was  introduced in the fifties in the context of the usual  games by Nash, see for example \cite{Na}.

\begin{definition}\label{def-Nash-measures}
We say that a pair  $(\bar \mu, \bar \nu)\in \stone \times \sttwo$ is a {\it Nash equilibria} for the game described  above if $\bar \mu \in BR_1(\bar \nu)$ and $\bar \nu \in BR_2(\bar \mu)$.
\end{definition}

If $(\bar \mu, \bar \nu)\in \stone \times \sttwo$ is a  Nash equilibrium point, we have that:

\medskip

 (a) player 1 can not get a better result by playing a different strategy $\tilde{\mu} \neq \bar{\mu}$ {\it if} player 2 uses strategy $\bar \nu$, and 
\smallskip 
 
 (b) %analogous interpretation for player 2. 
 player 2 can not get a better result by playing a different stratey $\tilde{\nu} \neq \bar{\nu}$ {\it if} player 1 uses strategy $\bar \mu$. 

\medskip
 
 So unilateral changes of strategy are not welcome, which makes $(\bar \mu, \bar \nu)$ a {\it stable} choice of strategies for both players.

Before passing to examples of ergodic games, let us remark that, in real life, Nash equilibrium is verified in many situations. For an example of Nash equilibrium in usual games, suppose a certain set of companies share the market of telephonic communications in a certain country. Each company has its prices for a certain set of services all of them offer to the costumers. Suppose company A analyse the market and conclude it can not get a better profit by changing its prices (supposing the other companies keep their prices unchanged). Then company A has no reason to modify its prices (its strategy). If all other companies are in the same situation, then the Nash equilibrium is attained, and is stable in the sense that no players will change their strategies {\it unilaterally}.

Now we consider  two examples of ergodic games, in the case  payoffs are given by \eqref{potenciais}:
\begin{itemize}
\item[a)] $A_2=-A_1$. This is called a zero-sum game, where 
the gains of one player are based on the losses of the other.
In this particular case we can use a minimax formulation to define the Nash equilibrium of the game. %\footnote{referência?}. 
Zero sum games are well-known examples of games, but are very restrictive.  For this reason, we will not consider the minimax formulation here.% , that can be only used in zero-sum games. %We prefer to use the more general 
\item[b)] $A_2=A_1$. In this case, one could think we have here a simple problem of ergodic optimization (see \cite{Jen}), and indeed there is some relation with ergodic optimization, in the sense that a choice $(\bar \mu,\bar \nu)$ of strategies that maximizes the common payoff is a {\it Nash equilibria}:
suppose $(\bar \mu, \bar \nu)$ maximizes 
\begin{equation}\label{payoffsiguais}
\int_{X\times Y} A_1 (x,y) d\mu(x) d\nu(y) =  \int_{X\times Y} A_2 (x,y) d\mu(x) d\nu(y)
\end{equation} 
over the set of mutual strategies $\stone \times \sttwo$ (which is a product space). Then, as a direct consequence of definition \ref{def-Nash-measures}, 
$(\bar \mu, \bar \nu)$ is a Nash equilibrium for the ergodic game.
 However, as we just consider here product measures, this is {\it not} an  ergodic optimization problem in the usual sense. Another feature of non-cooperative games which shows  that the searching for Nash equilibrium in the case $A_1=A_2$ can not be reduced to solving an ergodic optimization problem is the fact that, as there is no communication among players, they can not decide to use a common strategy. 
\end{itemize}

% {\bf Example:} Still in the case payoffs are given by \eqref{potenciais}, a simple example is given by the case $A_1=A_2$. 

%Podemos provar a existência do equilíbrio de Nash usando o teorema de ponto fixo de Kakutani (o mesmo usado em GMS10) enquanto a unicidade pode ser obtida via hipóteses adequadas de convexidade em $A_1$ e $A_2$.

%{\bf Example } BLABLABLA 
 
% \newpage
 
\section{Ergodic optimization and games}
 
In this section, we suppose
payoffs are given by \eqref{potenciais}. Although  stated for the case of two players, the following results can be easily generalized to any finite number of players. 
 
\begin{theorem}\label{existencia}
Suppose the payoffs are given by \eqref{potenciais}, where the potentials are Lipschitz-continuous functions. Then, there exists a Nash equilibrium. 	
\end{theorem} 
 
The proof of this theorem is an application of the Kakutani-Fan-Glicksberg Theorem
(see \cite{kakpaper}). We begin with some preliminary results that allow us to introduce this set-valued fixed point theorem. Remember that, if $M$ is a set, $\mathcal{P}(M)$ denotes the set of subsets of $M$. 

\begin{definition}
Let $M$ be a set. A set-valued function on $M$ is a function $\xi: M \to \mathcal{P}(M)$ such that $\xi(x)\neq \emptyset$ for all $x\in M$.
\end{definition}

%\newpage

\begin{definition}\label{KakutaniMap}
Let $E$ be a Hausdorff topological vector space, $K$ a non-empty compact and convex subset of $E$, and $\xi: K \to \mathcal{P}(K)$ a set valued function on $K$. Then $\xi$ is a Kakutani map if:
\begin{enumerate}
\item[a)]  $\forall x \in K$, $\xi(x)$ is a closed and convex non-empty subset of $\mathcal{P}(K)$.
\item[b)]  the set
$$ Gr(\xi)=\{(x,y)\;|\; x\in K, y \in \xi(x)\}$$
is closed in the product topology of $K\times K$. This means that if $(x_n,y_n) \to (x,y)$ and $y_n \in \xi(x_n)$ then $y \in \xi(x)$.
\end{enumerate}
\end{definition}

\begin{theorem}[Kakutani-Fan-Glicksberg]\label{KakutaniTheo}
Let $K$ be a non-empty compact and
convex subset of a convex Hausdorff topological vector space.
Let $\xi:K\to \mathcal{P}(K)$ be a Kakutani map. Then there exist a point $x^*\in K$ such that $$x^* \in \xi(x^*).$$
\end{theorem}

Before we are able to prove our first existence result, we need to prove a technical lemma. First remember the Wasserstein-1 metric over $\mathcal{M}(X)$,  given by
\begin{equation}\label{WassM}
W^1(\mu_1,\mu_2)=\sup_{Lip(\varphi)\leq 1} \int \varphi d(\mu_1-\mu_2),
\end{equation}
where
\begin{equation*}\label{lipnorm}
Lip(\varphi)=\sup_{x \neq x'}\frac{|\varphi(x)-\varphi(x')|}{d_X(x,x')}
\end{equation*}
 is the Lipschitz norm on the set of Lipschitz functions from $X$ to $\mathbb{R}$.
 This  form of the Wasserstein-1 metric is a well-known consequence of Kantorovich-Rubinstein duality theorem (see \cite{Vi}). 
 %page 207).  
 As we are in compact spaces, we know $W^1$ metrizes the weak-* topology. % (see also \cite{Vi} page 212). 
We also have that, for any Lipschitz map $\psi:X \to \mathbb{R}$,
\begin{equation}\label{desig-w1}
\int \psi d(\mu_1 - \mu_2) \leq Lip(\psi) W^1(\mu_1,\mu_2).
\end{equation}

\begin{lemma}\label{technical-for-kakutani}
Suppose $\nu_n$ and $\hat{\mu}_n$ are sequences of probabilities on $Y$ and $X$, respectively, and that $\nu_n \to \nu$ and $\hat{\mu_n} \to \hat{\mu}$.
Then, if $A_1$ is Lipschitz-continuous, we have 

$$\int A_1(x,y) d\hat{\mu}_n(x) d\nu_n(y) \to \int A_1(x,y) d\hat{\mu}(x) d\nu(y)\,.$$

\end{lemma}

\noindent
{\bf Proof of lemma \ref{technical-for-kakutani}:}
We use Fubini Theorem:
$$\left|\int A_1(x,y) d\hat{\mu}_n(x) d\nu_n(y) - \int A_1(x,y) d\hat{\mu}(x) d\nu(y) \right| \leq $$
$$\leq \left|\int \int A_1(x,y)  d\nu_n(y) d\hat{\mu}_n(x) - 
\int \int A_1(x,y)  d\nu(y) d\hat{\mu}_n(x)\right| +$$
$$+  \left|\int \int A_1(x,y)  d\nu(y) d\hat{\mu}_n(x)
- \int \int A_1(x,y)  d\nu(y) d\hat{\mu}(x) \right|. $$

Now, compactness implies that $x \mapsto \int A_1(x,y)  d\nu(y)$ is a continuous function, and this means that the second term above converges to zero when $n \to \infty$. The first term also goes to zero because
it is bounded by
$$\int \left| \int A_1(x,y)  d\nu_n(y)  - \int A_1(x,y)  d\nu(y) \right| d\hat{\mu}_n(x), $$
and using \eqref{desig-w1} we have, $\forall x$, 
$$ 
\left| \int A_1(x,y)  d\nu_n(y)  - \int A_1(x,y)  d\nu(y) \right|  \leq Lip(A_1) W^1(\nu_n,\nu). 
$$ 
\cqd 
%$$\left|\int \int A_1(x,y)  d\nu_n(y) d\hat{\mu}_n(x) - \int \int A_1(x,y)  d\nu(y) d\hat{\mu}_n(x)+\right| = $$
%$$ = \left|\int \int A_1(x,y)  d\nu_n(y) d\hat{\mu}_n(x) - \int \int A_1(x,y)  d\nu(y) d\hat{\mu}_n(x)+\right|

\medskip

\noindent{\bf Proof of theorem \ref{existencia}:}
 Let $K=\stone \times \sttwo$. Using lemma \ref{primeiro}, we know that $K$ is a non-empty, compact and convex subset of the Hausdorff topological vector space given by the product of the spaces of  finite signed %regular Borel 
 measures on $X$ and $Y$. 
%$\mathcal{M}(X) \times \mathcal{M}(Y)$% endowed with the product topology.
%the product of the sets of probability measures on $X$ and $Y$, with the product topology.

%%%We know that the set of invariant probability measures for a dynamics defined on a compact metric space $C$ is a non-empty compact and convex subset of the set of probability measures on the borelian sets of $C$. Therefore, $K$, being the product of two of such sets, is also a non-empty, compact and convex subset of the Hausdorff topological vector space given by the product of the sets of probability measures on $X$ and $Y$, with the product topology.

Now let $\xi=BR:K \to \mathcal{P}(K)$ be the set-valued function given by:

$$BR(\mu,\nu)= BR_1(\nu) \times BR_2(\mu).$$

{\it Claim:} $BR$ is a Kakutani map (see definition \ref{KakutaniMap}).

\medskip
{\it Proof of the Claim:}
We begin by noting that the condition (a) of definition \ref{KakutaniMap} is a direct consequence of lemma \ref{BR_properties}.

Now we address condition (b) of definition \ref{KakutaniMap}. 
Suppose $\Big((\mu,\nu),(\hat{\mu},\hat{\nu})\Big)$ belongs to the closure of $Gr(BR) \subset K \times K$. We need to prove that  $(\hat{\mu},\hat{\nu})\in BR(\mu,\nu)$. 
%$(\mu,\nu) \in \overline{Gr(\xi)}$. 
To do that, we know 
$K=\stone \times \sttwo$ is metrizable and so is $K \times K$. Therefore,  there exists a sequence 
$\Big((\mu_n,\nu_n),(\hat{\mu}_n,\hat{\nu}_n)\Big)$ converging to 
$\Big((\mu,\nu),(\hat{\mu},\hat{\nu})\Big)$, such that 
$(\hat{\mu}_n,\hat{\nu}_n)\in BR(\mu_n,\nu_n)$.
Now 
$\hat{\mu}_n\in BR_1(\nu_n)$ for all $n\in \nat$. This means that for any other strategy $\tilde \mu$ for player 1, we have
$$\varphi_1(\hat{\mu}_n,\nu_n) \geq \varphi_1(\tilde \mu,\nu_n),$$ 
which means
$$\int A_1(x,y)d\hat \mu_n(x)d\nu_n(y) \geq \int A_1(x,y)d\tilde \mu(x)d\nu_n(y).$$

If we take the limit when $n \to \infty$, using lemma \ref{technical-for-kakutani}, we get

$$\int A_1(x,y)d\hat\mu(x)d\nu(y) \geq \int A_1(x,y)d\tilde \mu(x)d\nu(y)\,\;\forall \, \tilde \mu \in \stone,$$
%As $\tilde \mu$ was arbitrary, 
which means that 
$\varphi_1(\hat{\mu},\nu) \geq \varphi_1(\tilde \mu,\nu)\,\;\forall \, \tilde \mu \in \stone,$
and this implies $\hat \mu \in BR_1(\nu)$.

In an analogous way we can prove that  $\hat \nu \in BR_2(\mu)$, and these two facts imply that $(\hat\mu,\hat\nu)\in BR(\mu,\nu)$.
 This ends the proof of the {\it Claim}.

\medskip

Now Kakutani Fixed Point (theorem \ref{KakutaniTheo}) implies the existence of  $(\bar \mu, \bar \nu)$ satisfying $(\bar \mu, \bar \nu)\in 
BR(\bar \mu, \bar \nu)$ . Therefore, 
$$
\begin{cases}
\bar \mu \in BR_1(\bar \nu), \\ \bar \nu \in BR_2(\bar \mu),
\end{cases}
$$
and this means that $(\bar \mu, \bar \nu)$ is a Nash equilibrium point.
\cqd
%{\it End of Proof of Theorem \ref{existencia}.}

%\end{document}

\section{ Thermodynamic formalism and games}\label{sec:entropia}

In this section, we suppose the payoffs are given by \eqref{payoff-entropia}, i.e.

\begin{equation*}
\begin{cases}
\varphi_1(\mu,\nu)=\int_{X \times Y} A_1(x,y) d\mu(x) d\nu(y)+h(\mu),
\\
\varphi_2(\mu,\nu)=\int_{X \times Y} A_2(x,y) d\mu(x) d\nu(y)+h(\nu),
\end{cases}
\end{equation*}
where $h(\mu)$ and $h(\nu)$ are the metric entropy of the measures $\mu$ and $\nu$. 
We also suppose the potentials $A_1$ and $A_2$ are Lipschitz maps and  satisfy the condition
\begin{equation}\label{condition}
|A_i(x,y)-A_i(x',y)-A_i(x,y')+A_i(x',y')|<C d_X(x,x') d_Y(y,y')
\;\; \forall \; i=1,2,\end{equation}
where $C>0$. Condition \eqref{condition} holds, for example, if $X=Y=[0,1]$ and $A_1$ and $A_2$ are $C^2$ maps.
%We remark that condition \eqref{condition} holds in numerous situations. For instance, if $X=Y=[0,1]$ and $A_1$ and $A_2$ are $C^2$ maps.

Finally, we suppose $X$ and $Y$ are both given by the Bernoulli space of sequences in $d$ symbols, 
i.e. $$X=Y= \mathcal{B}\equiv\{1,2,3\hdots,d\}^{\nat},$$ with the usual metric (see \cite{PP})
and the maps $T$ and $S$ are given by the shift map on $\mathcal{B}$. 
%the maps $T:X\rar X$ and $S:Y \rar Y$ are given by the shift map and expanding continuous maps of the circle, i.e. $X=Y=S^1$,  and there exist $\sigma_i>1$ and $\delta_i>0$, $i=1,2$, such that  $d(T(x),T(y))>\sigma_1 d(x,y)$  if $d(x,y)<\delta_1$, and $d(S(x),S(y))>\sigma_2 d(x,y)$   if $d(x,y)<\delta_2$.%, (same applies to $S$).

%{\footnote Definir a metrica? acho que sim...}

Remark: As an alternative to this hypothesis, we could require the maps $T$ and $S$ to be expanding continuous maps of the circle, or more generally, topologically mixing expanding continuous maps 
on compact metric sets. All the results of this section hold under any one of these other hypothesis.
In fact all we need is to be able to define the Ruelle operator and equilibrium measures, and also have  uniqueness of equilibrium measures  (see lemma \ref{TFG}, and references  \cite{Ka,PP}).
%and the continuous dependence of the equilibrium measure on the potential

The main result of this section is

\begin{theorem}\label{existence-FT}
Let $T$ and $S$ be the shift map on the Bernoulli space of $d$ symbols. %expanding continuous maps of the circle. 
Suppose the potentials $A_1$ and $A_2$ are Lipschitz, and satisfy \eqref{condition}. Then there exist a Nash equilibrium point for the game where payoffs are given by \eqref{payoff-entropia}.
\end{theorem}

% \subsection{Uniqueness}
%Here is the novelty, over the case of subsection \ref{uniq-EO}. 

%It is possible to give a proof of this result with the use of Kakutani fixed point. Instead, we prefer to give a proof based on arguments of thermodynamic formalism. 
 We begin with some technical lemmas.

\medskip

%\subsection{A technical result}

%%
%%%Let $\mu_1$ and $\mu_2$ be two probability measures on $X$. We have
%%
%%%$$ \int \varphi d(\mu_1-\mu_2) \equiv  \int \varphi d\mu_1-\int \varphi d\mu_2.$$
%%%Note that $\mu_1-\mu_2$ is a signed measure.  %Note that %the following is FALSE:
%%%$$\int \varphi d(\mu_1-\mu_2) \leq \int |\varphi| d(\mu_1-\mu_2) \;\;\mbox{CAN BE FALSE}.$$

%\footnote{repetido !!!}
%Let us remember the Wasserstein-1 metric over $\mathcal{M}(X)$, that has a particular form given by
%$$W^1(\mu_1,\mu_2)=\sup_{Lip(\varphi)\leq 1} \int \varphi d(\mu_1-\mu_2).$$
% This particular form of the Wasserstein-1 metric is a well known consequence of 
%Kantorovich-Rubinstein duality theorem (see \cite{Vi} page 207).
% As we are in compact spaces, we know $W^1$ metrizes the weak-* topology
% (see also \cite{Vi} page 212).

Remember that $\mathcal{M}(X)$ and $\mathcal{M}(Y)$ are metric spaces with the Wasserstein-1 metric \eqref{WassM} and that the topology generated by this metric is equivalent to the weak-* topology.
Let us also denote  by $\lip$ the set of Lipschitz maps from $X$ to $\mathbb{R}$. 
We define in $\lip$ the norm $\| \cdot \|:\lip \to \re$ given by
$$\|\psi\|=\|\psi\|_0+Lip(\psi),$$ 
where $\|\psi\|_0=\sup_{x\in X}|\psi(x)|$ is the usual $C^0$ norm
and %\begin{equation}\label{lipnorm}
$Lip(\psi)=\sup_{x \neq x'}\frac{|\psi(x)-\psi(x')|}{d_X(x,x')}$
%\end{equation}
 is the Lipschitz norm. 
It is a well-known fact that $\|\psi\|$ is  a norm that makes $\lip$ a Banach space (see \cite{PP}). In an analogous way, we define $\lipy$ and its norm. 
%, with the same (kind of) norm.

\begin{lemma}\label{l1} 
Suppose $A_1:X\times Y \to \re$ is Lipschitz and also satisfies  condition \eqref{condition}.
For any $\nu\in \mathcal{M}(Y)$, 
the  function  $\psi_{\nu}:X \to \mathbb{R}$ given by 
\begin{equation}\label{psinu}
\psi_{\nu}(x)=\int_Y A_1(x,y)d\nu(y)
\end{equation}
is Lipschitz. Also, 
the map that sends $\nu\in \sttwo$ to $\psi_{\nu}\in \lip$, given by \eqref{psinu}, is a Lipschitz continuous map: we have 
\begin{equation}\label{lip-psi}
 \|\psi_{\nu}-\psi_{\nu'}\|\leq \left(C+ Lip(A_1)\right) W^1(\nu,\nu') .
\end{equation}
\end{lemma}
\medskip

{\bf Proof of lemma \ref{l1}}: 
The map \eqref{psinu} is clearly Lipschitz. 
In order to prove \eqref{lip-psi}, we 
claim that  for any pair of points $x$ and $x'$ in $X$, the function $\xi_{x,x'}:Y \to \mathbb{R}$ given by 
 $$\xi_{x,x'}(y)=A_1(x,y)-A_1(x',y)$$ is Lipschitz and $Lip(\xi_{x,x'})=C d_X(x,x')$.

\medskip

This claim is a direct consequence of the hypothesis \eqref{condition}: for any $y$ and $y'$, we have
\begin{align*}
|\xi_{x,x'}(y)-\xi_{x,x'}(y')| & = |A_1(x,y)-A_1(x',y)-A_1(x,y')+A_1(x',y')|\\
& < C d_X(x,x') d_Y(y,y').\end{align*}

Now we prove \eqref{lip-psi}: we begin by considering the $C^0$ norm: note that 
\begin{align*}
|\psi_{\nu}(x)-\psi_{\nu'}(x)| & =  \left|\int_Y A_1(x,y)  d\nu(y)-\int_Y A_1(x,y)  d\nu'(y)\right|\\
&= \left|\int_Y A_1(x,y)  d(\nu-\nu')(y)\right|\\ 
&\leq Lip(A_1) W^1(\nu,\nu'),
\end{align*}
where we used \eqref{desig-w1}. 
This proves that
$$\|\psi_{\nu} - \psi_{\nu'}\|_0 =\sup_{x\in X}
|\psi_{\nu}(x)-\psi_{ \nu'}(x)|\leq 
Lip(A_1) W^1(\nu,\nu') .$$

To deal with the Lipschitz norm, we note that
\begin{align*}
|(\psi_{\nu}-\psi_{\nu'})(x)-(\psi_{\nu}-\psi_{\nu'})(x')| & = \Bigg|\int_Y(A_1(x,y)-A_1(x',y))d(\nu-\nu')(y)\Bigg| \\
&   \leq Lip(\xi_{x,x'})W^1(\nu,\nu') 
\leq C d_X(x,x')W^1(\nu,\nu') ,\end{align*}
where we used \eqref{desig-w1} again, and the claim.
This implies that
$$Lip(\psi_{\nu} - \psi_{\nu'}) = \sup_{x\neq x'}
\frac{\|(\psi_{\nu}-\psi_{ \nu'})(x)-(\psi_{\nu}-\psi_{ \nu'})(x')\|}{d_X(x,x')}\leq   C  W^1(\nu,\nu'),$$
which finishes the proof of the lemma.
\cqd

The next lemma is  analogous to lemma \ref{l1}:

\begin{lemma}\label{l2}
Suppose $A_2:X\times Y \to \re$ is Lipschitz and satisfies condition \eqref{condition}.

For any $\mu\in \mathcal{M}(X)$, 
the  function  $\psi_{\mu}:Y \to \mathbb{R}$ given by 
\begin{equation}\label{psinu-2}
\psi_{\mu}(y)=\int_X A_2(x,y)d\mu(x)
\end{equation}
is Lipschitz. Also, 
the map that sends $\mu\in \stone$ to $\psi_{\mu}\in \lipy$, given by \eqref{psinu-2}, is a Lipschitz continuous map: we have 
\begin{equation}\label{lip-psi-2}
 \|\psi_{\mu}-\psi_{\mu'}\|\leq \left(C+ Lip(A_2)\right) W^1(\mu,\mu') .
\end{equation}
\end{lemma}

\medskip

Before passing to the last lemma, we need 
to recall some facts about  equilibrium measures and the  Ruelle-Perron-Frobenius operator (see \cite{Ka, OV, PP}). % This facts will be also used in the next section.

The Ruelle-Perron-Frobenius operator associated to a Lipschitz {\it potential}  $A:X \to \mathbb{R}$ is the operator on $\lip$ that associates to any $\varphi \in \lip$ the function $L_A(\varphi) \in \lip$ given by 

$$L_A(\varphi)(x)=\sum_{T(z)=x}e^{A(z)}\varphi(z).$$

Remember that if $T$ is the shift map on the Bernoulli space of symbols, or if it is an expanding continuous map of the circle, then it has finite degree, i.e. the set $\{z\in S^1 \| T(z)=x\}$ is finite and has the same number of elements for any $x$
(this result also holds when $T:X\rar X$ is topologically mixing expanding continuous map on a compact metric space $X$. See, for instance, \cite{ Ka, OV}.)%AS already said, all the result

%(Even under the more general hypothesis where $T:X\rar X$ is topologically mixing expanding continuous map on a compact metric space $X$, any point has a finite number of inverse images, and the all the results here still hold. See, for instance, \cite{ Ka, OV}.)

%Under our hypothesis ($T:S^1 \to S^1$ an expanding continuous map of finite degree)
%%%($T:X\rar X$ a topologically mixing expanding continuous map on the compact metric space $X$, and $A$ Lipschitz continuous)
 %the RPF operator is well defined, i.e., the sum defining it is a convergent sum, because it is a finite sum. % - see for example \cite{}). ACHO QUE DA PARA CITAR O MARCELO VIANA / KRERLEY - MAS ESTA EM PORTUGUES...%CUIDADO O KATOK ATE FALA ISSO MAS DE MANEIRA MAIS COMPLICADA...

A function $A\in \lip$ is called a {\it normalized potential} if $L_A(1)=1$. The RPF operator has a maximal eigenvalue $\lambda_A>0$ associated to an eigenfunction $\varphi_A$, which is simple and positive (simple means the eigenspace has dimension $1$). %all other eigenfunctions are multiples of $\psi_A$). %\footnote{normalization in the case of the eigenfunction means BLABLABALBA}
If $A$ is non-normalized then it can be normalized by considering
$$\bar A = A + \log \varphi_A-\log \varphi_A \circ T - \log \lambda_A.$$

The dual RPF operator, denoted by $L_A^*$, acts on probability measures on $X$, and is defined by $$\int \psi d L_A^*(\mu) = \int L_A(\psi) d\mu.$$

If $\bar A$ is the normalized potential associated to $A$, the dual operator $L_{\bar A}^*$ has a unique fixed probability $\mu_A$ (i.e. a probability $\mu_A$ that satisfies $L_{\bar A}^*(\mu_A)=\mu_A$), called the Gibbs state associated to $A$, which is invariant for $T$ (and is also exact and ergodic), and satisfies $$ \int A d \mu + h(\mu) \leq \int A d{\mu_A} + h(\mu_A)\;,\; \forall \;\; \mu \in \stone.$$
Moreover, $\mu_A$ is the unique invariant measure to attain the maximum above. A measure that attains the maximum above is called the equilibrium measure for $A$, is unique, exact and ergodic and is the Gibbs measure $\mu_A$. The last result is also known as the {\it variational principle for pressure}
(see also \cite{Ka} chapter 20.3 or \cite{OV}).

A final and very important fact about equilibrium measures (see \cite{Ma} and also \cite{ KLS, Varandas} for more recent results on more general settings) is the fact that
the  function that associates to $A \in \lip$ the equilibrium measure $\mu_A \in \stone$ is a continuous map (in fact, such function is even analytic - but here we will only need its  continuity).
%Some final and very important facts about the RPF operator (see \cite{PP}. See also \cite{Varandas,SSS} for some more recent works in more general settings): 
%\begin{itemize}
%\item the  functions that associate to $A \in \lip$ the maximal eigenvalue $\lambda_A \in \mathbb{R}$ and the associated normalized eigenfunction $\varphi_A \in \lip$ are analytical functions.
%%%% defined in the Banach space $\lip$.
%\item the  function that associates to $A \in \lip$ the equilibrium measure $\mu_A \in \stone$ is a continuous map. 
%\end{itemize}

Now, we are able to prove the next lemma:

\begin{lemma}\label{TFG} 
%Let $T$ and $S$ be continuous topologically mixing expanding continuous maps, 
%Let $T$ and $S$ be expanding continuous maps of the circle,
Let $T$ and $S$ be the shift map on the Bernoulli space of $d$ symbols, 
and suppose that $A_1$ and $A_2$ are Lipschitz continuous maps that satisfy  \eqref{condition}.  
%\footnote{como vou usar a métrica wasserstein-1, preciso Lipschitz. Acho que Holder nao vai ser suficiente.}
 Then:
\begin{enumerate}
\item[a)] $BR_1(\nu)$ and $BR_2(\mu)$ both have only one measure (are singleton sets). 
\item[b)] $BR_1:\mathcal{M}_S(Y)\to \mathcal{M}_T(X)$ is a continuous function (in the weak-* topology).
\item[c)] $BR_2:\mathcal{M}_T(X)\to \mathcal{M}_S(Y)$ is a continuous function (in the weak-* topology).
\end{enumerate} 
\end{lemma}

\noindent{\bf Proof of lemma \ref{TFG}:} (a) Let us fix $\nu \in \mathcal{M}_S(Y)$. We have 
\begin{equation}\label{BR1}
BR_1(\nu)= \verb"argmax" \left\{ \mu \mapsto 
\int_{X} \Big(\int_Y A_1(x,y)  d\nu(y)\Big) d\mu(x) +h(\mu)
\right\}.
\end{equation}

Now, if we remember that 
\begin{equation*}%\label{psinu}
\psi_{\nu}(x)= \int_Y A_1(x,y)  d\nu(y),
\end{equation*}
using lemma \ref{l1}, we have that  $\psi_{\nu}$ is a Lipschitz potential (defined on $X$) and any measure on $BR_1(\nu)$ is an equilibrium measure for $(T,\psi_{\nu})$. If we use the uniqueness of equilibrium measures, we get that $BR_1(\nu)$ is a singleton set. Analogously, using lemma \ref{l2} we prove that $BR_2(\mu)$ has only one measure, which finishes the proof of the first item of the lemma.

(b) We just need to use lemma \ref{l1}, that proves that the map 
$\nu \in \mathcal{M}_S(Y)\mapsto \psi_{\nu}\in \lip$ is a Lipschitz map,  and the fact that the equilibrium measure depends continuously 
 on the potential $\psi_{\nu}$. % (see \cite{Ma}).
%In fact we have that the map $BR_1$ is a Lipschitz map.

%$$\psi_{\nu_n}(x)= \int_Y A_1(x,y)  d\nu_n(y) \to \psi_{\nu}(x)= \int_Y A_1(x,y)  d\nu(y).$$

(c) It is analogous to (b), using lemma \ref{l2}.%We just have to use the same idea of item (b).
\cqd

\vspace{0.3cm}

%\subsection{Existence }

\noindent{\bf Proof of theorem \ref{existence-FT}:}
%Here we could have used use Kakutani Theorem, but we prefer to give a direct proof based on lemma \ref{TFG}.

Define 
$BR: \mathcal{M}_T(X)\times \mathcal{M}_S(Y) \to  \mathcal{M}_T(X)\times \mathcal{M}_S(Y)$
as
\begin{equation*}%\label{BRmap}
BR(\mu,\nu)=BR_1(\nu) \times BR_2(\mu).
\end{equation*} 

As a consequence of lemma \ref{TFG}, this is a well defined map. Moreover, it is continuous, because both functions $BR_1$ and $BR_2$ are continuous. Now we can use Tychonoff-Schauder fixed point theorem to get a fixed point, which is a Nash equilibrium. \cqd

%This ends the proof of the Theorem.

%\subsection{Uniqueness}\label{uniq-TF}

%\newpage

\section{Uniqueness}\label{uniq-EO}%\footnote{olhar !!}

%Uniqueness of Nash equilibrium is a desirable property of any game. In the case there is only one Nash equilibrium both players can predict the other behaviour and then we can expect that the 
%TEXTO NOVO:

Now, we make some considerations about uniqueness of Nash equilibrium. We remark that, even if we have uniqueness for maps $BR_1$ and $BR_2$ (i.e. even if the best response of player one, given any strategy of player two, is unique - and vice versa) this does not imply uniqueness of Nash Equilibrium. % This is due to the fact that players do not communicate, and so player one a priori can not predict or negotiate the behaviour of player two. 
We can have more than one Nash equilibrium point even in this case: it is possible that $(\bar \mu, \bar \nu)$ and $(\tilde \mu, \tilde \nu)$ are two Nash equilibrium points, with a unique  best response of both players in both situations. %We now present two examples. The second one presents lack of uniqueness with unique best response of both players.
%%, That is exactly what happens in the next example:

%Suppose we are in the case $A_1=A_2$: suppose there are two pairs $(\bar \mu, \bar \nu)$ and $(\tilde \mu, \tilde \nu)$ that maximize \eqref{payoffsiguais}. Then there is (at least) two Nash equilibria for this ergodic game.

%Now we pass to a second example where $A_1 \neq A_2$.

Let us analyse an example. Let $X=Y=\mathcal{B}$ be the Bernoulli space, now on two symbols $0$ and $1$.
An element of $\mathcal{B}$ is denoted by $x$, where % and is a sequence of $0$ and $1$, 
$x=(x_1x_2x_3x_4\hdots)$, with $x_i \in \{0,1\}$ for any $i\geq 1$. 
%The distance map is given by  $d(x,y)= \sum_{i \geq 1} \frac{|x_1-y_1|}{2^i}$.
We denote by $\bar 0$ the element $\bar 0 = (0000\hdots)$ and 
$\bar 1$ the element $\bar 1 = (1111\hdots)$.
The Dirac measure $\delta_{\bar 0}$ in $\bar 0$ is an invariant measure for the shift, as well as the 
Dirac measure $\delta_{\bar 1}$ in $\bar 1$.

Now let $A_1$ and $A_2$ be any two Lipschitz maps such that 
$$
\begin{cases}
A_1(\bar 0 ,\bar 0)=2 ,\\
A_1(x ,\bar 0) <2 \mbox{ if }x \neq \bar 0 ,\\
A_2(\bar 0 ,\bar 0)=3 ,\\
A_2(\bar 0,y)<3  \mbox{ if }y \neq  \bar 0,
\end{cases}
$$
and
$$
\begin{cases}
A_1(\bar 1 ,\bar 1)=3 ,\\
A_1(x ,\bar 1) <3 \mbox{ if }x \neq \bar 1 ,\\
A_2(\bar 1 ,\bar 1)=2, \\
A_2(\bar 1,y)<2  \mbox{ if }y \neq \bar 1.
\end{cases}
$$
Suppose the payoff functions are given by \eqref{potenciais}.
% A_1 é o ganho do homem, 0=ballet, 1=football. 
 We see that the two pairs $\delta_{\bar 0} \times \delta_{\bar 0}$ 
and  $\delta_{\bar 1} \times \delta_{\bar 1}$ 
% the Dirac measures on both pairs $(\bar 0 ,\bar 0)$ and $(\bar 1 ,\bar 1)$ 
 are Nash equilibrium, and also that both maps $BR_1$ and $BR_2$ are single-valued in both pairs, i.e, best response of both players is unique in both situations:
$$ 
\begin{cases}
BR_1(\delta_{\bar 0}) = \delta_{\bar 0},\\
BR_1(\delta_{\bar 1}) = \delta_{\bar 1},\\
BR_2(\delta_{\bar 0}) = \delta_{\bar 0},\\
BR_2(\delta_{\bar 1}) = \delta_{\bar 1}.
\end{cases}
$$

\section{Ergodic Transport and a cooperative game }

This section is independent of the results of the preceding sections. %\footnote{dá uma ajeitada rapida - e tenta falar um pouco sobre jogos nao cooperativos.}
%\bigskip
Here, we will introduce a cooperative game. We will use tools from ergodic transport (see \cite{LM}).
As we said in the introduction, we will not elaborate on this example. We just present it and plan to analyse it further in a future work.

%Remember we have two dynamics $T:X\to X$ and $S:Y \to Y$.
Suppose the  payoffs $\varphi_1:X \times Y \to \re$
and $\varphi_2:X \times Y \to \re$ are given by \eqref{potenciais}. 
Suppose, however, that player´s 1 payoff does not depend on the strategy of player 2. This holds when $A_1(x,y)=A_1(x)$ is a function of the first variable only. 
In this case player 1 chooses his optimal strategy by solving a simple ergodic optimization problem: Let $\mu$ be an invariant measure that maximizes 
$$m \in \mathcal{M}_{T}(X) \to \int_X A_1(x) dm(x).$$ 

Now we consider player 2. Suppose $A_2(x,y)$ does depend on both variables,
which means player 2 depends on player's 1 choices.
On his side, player 1 accept to collaborate with player 2 and both search for a common strategy $\pi \in \mathcal{M}_{T,S}(X \times Y)$ (which denotes the set of plans whose x-projection is $T-$invariant  and y-projection is $S-$invariant)
provided his expected payoff $\int A_1 d\mu$ does not gets smaller. This means that $\pi$ has to have x-projection given by $\mu$. In other words,  $\pi$ needs to belong to  $\in \mathcal{M}_{\mu,S}(X \times Y)$, which is the set of plans whose x-projection is $\mu$ and y-projection is invariant for $S$ (see \cite{LM}).

So what they need to do is to solve an ergodic transport problem:  they search for the common strategy $\pi^*$ that maximizes  

$$\pi \in \mathcal{M}_{\mu,S}(X \times Y) \to \int_{X \times Y} A_2(x,y)d \pi(x,y).$$

Here, we have a cooperative game because the common strategy is not a product measure as in the preceding sections. % This means they do not play independently.

Ergodic transport was introduced recently, see for example \cite{LM,LMMS2} and among other applications, we can find the invariant measure that minimizes the Wasserstein-2 distance (see \cite{Vi}) from a given measure which is not invariant: we just have to use as $A_2$ the cost function given by the square of the distance, and search for the measure that minimizes (instead of maximize) the expression above. In this case, $\mu$ will play the role of any given measure, not necessarily invariant.

%\footnote{falar um pouquinho mais do LM.}
%\footnote{a versao com entropia não seria uma aplicação do LMMS? Ainda, se citar o LMMS não precisa se preocupar em maximizar ao inves de minimizar. MAS lembra que em LMMS um dos espaços é finito.}

{\sc Rafael Rigão Souza\\
Departamento de Matematica\\
Universidade Federal do Rio Grande do Sul - UFRGS\\
Avenida Bento Gonçalves, 9500\\
91509-900 - Porto Alegre RS\\
Brazil}\\
{\it email rafars@mat.ufrgs.br}

\end{document}